\documentclass[12pt]{article} 
\usepackage{amssymb,amsmath}
\usepackage{amsthm}
\input{psfig.sty}

\newtheorem{theorem}{Theorem}
\newtheorem{proposition}[theorem]{Proposition}

\newtheorem{conjecture}[theorem]{Conjecture}

\theoremstyle{definition}
 
\newtheorem{example}{Example}

\theoremstyle{remark}

\begin{document}
\newcommand{\beq}{\begin{equation}} 
\newcommand{\eeq}{\end{equation}}
\newcommand{\beas}{\begin{eqnarray*}} 
\newcommand{\eeas}{\end{eqnarray*}}
\newcommand{\zz}{\mathbb{Z}}
\newcommand{\pp}{\mathbb{P}} 
\newcommand{\nn}{\mathbb{N}}
\newcommand{\rr}{\mathbb{R}}
\newcommand{\bm}[1]{{\mbox{\boldmath $#1$}}}
\newcommand{\bmp}{\bm{p}}
\newcommand{\bmq}{\bm{q}}
\newcommand{\con}{\mathrm{Comp}(n)}
\newcommand{\sn}{\mathfrak{S}_n} 
\newcommand{\sk}{\mathfrak{S}_k} 
\newcommand{\snm}{\mathfrak{S}_n^{(m)}} 
\newcommand{\skm}{\mathfrak{S}_k^{(m)}} 
\newcommand{\fs}{\mathfrak{S}}
\newcommand{\st}{\,:\,} 
\newcommand{\as}{\mathrm{as}}
\newcommand{\is}{\mathrm{is}}
\newcommand{\lgn}{\mathrm{len}}
\newcommand{\dis}{\displaystyle}
\newcommand{\ptq}{p\times q}

\thispagestyle{empty}

\vskip 20pt
\begin{center}
{\bf A Conjectured Combinatorial Interpretation of the Normalized
  Irreducible Character Values of the Symmetric Group}\\ 
\vskip 15pt
{\bf Richard P. Stanley}\\
{\it Department of Mathematics, Massachusetts Institute of
Technology}\\
{\it Cambridge, MA 02139, USA}\\
{\texttt{rstan@math.mit.edu}}\\[.2in]
{\bf\small version of 21 June 2006}\\
\end{center}
The irreducible characters $\chi^\lambda$ of the symmetric group $\sn$
are indexed by partitions $\lambda$ of $n$ (denoted $\lambda\vdash n$
or $|\lambda|=n$), as discussed e.g.\ in \cite[{\S}1.7]{macd} or
\cite[{\S}7.18]{ec2}. If $w\in \sn$ has cycle type $\nu\vdash n$ then
we write $\chi^\lambda(\nu)$ for $\chi^\lambda(w)$. 

Let $\mu$ be a partition of $k\leq n$, and let $(\mu,1^{n-k})$ be the
partition obtained by adding $n-k$ 1's to $\mu$. Thus $(\mu,1^{n-k})
\vdash n$. Regarding $k$ as given, define the \emph{normalized
character} $\widehat{\chi}^\lambda(\mu,1^{n-k})$ by
  $$ \widehat{\chi}^\lambda(\mu,1^{n-k}) = \frac{(n)_k
    \chi^\lambda(\mu,1^{n-k})}{\chi^\lambda(1^n)}, $$
where $\chi^\lambda(1^n)$ denotes the dimension of the character
$\chi^\lambda$ and $(n)_k=n(n-1)\cdots (n-k+1)$. Thus
\cite[(7.6)(ii)]{macd}\cite[p.\ 349]{ec2} 
$\chi^\lambda(1^n)$ is the number $f^\lambda$ of standard Young
tableaux of shape $\lambda$. 

\begin{figure}
 \centerline{\psfig{figure=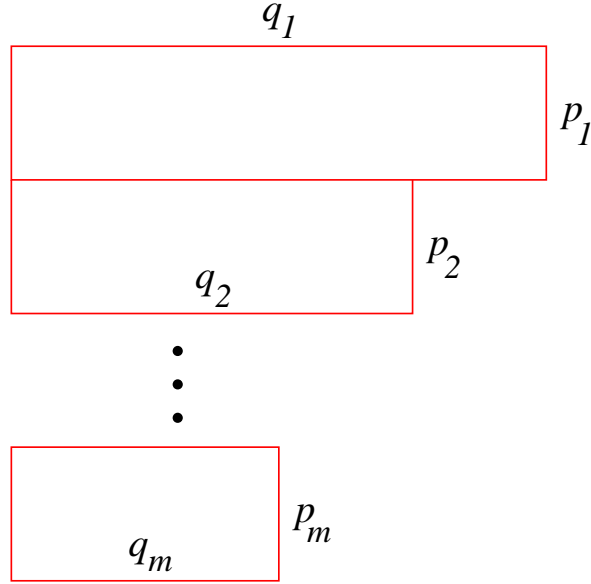}}
\caption{A union of $m$ rectangles}
\label{fig:nrec}
\end{figure}

Suppose that (the diagram of) the partition $\lambda$ is a union of $m$
rectangles of sizes $p_i\times q_i$, where $q_1\geq q_2 \geq \cdots
\geq q_m$, as shown in Figure~\ref{fig:nrec}. The following result was
proved in \cite[Prop.~1]{rc} for $\mu=(k)$ and attributed to
J. Katriel (private communication) for arbitrary $\mu$. 

\begin{proposition} \label{prop:f}
Let $\lambda$ be the shape in Figure~\ref{fig:nrec}, and fix $k\geq
1$. Let $\mu\vdash k$. Set $n=|\lambda|$ and 
  $$ F_\mu(\bmp;\bmq)=F_\mu(p_1,\dots,p_m; q_1,\dots, q_m) =
  \widehat{\chi}^\lambda(\mu,1^{n-k}). $$
Then $F_\mu(\bmp;\bmq)$ is a polynomial function of
the $p_i$'s and $q_i$'s with integer coefficients, satisfying
   $$ (-1)^k  F_\mu(1,\dots,1;-1,\dots,-1) = (k+m-1)_k. $$
\end{proposition}

\textsc{Note.} When $\mu=(k)$, the partition with a single part $k$,
we write $F_k$ for $F_{(k)}$. A formula was given in \cite[(9)]{rc}
for $F_k(\bmp;\bmq)$, viz., 
  $$ F_k(\bmp;\bmq) = -\frac 1k [x^{-1}]
   \frac{(x)_k\displaystyle \prod_{i=1}^m (x-(q_i+p_i+p_{i+1}+
    \cdots+p_m))_k}{\displaystyle \prod_{i=1}^m (x-(q_i+p_{i+1}
    +p_{i+2}+\cdots+p_m))_k}, $$
where $[x^{-1}]f(x)$ denotes the coefficient of $x^{-1}$ in the
expansion of $f(x)$ in \emph{descending} powers of $x$ (i.e., as a
Taylor series at $x=\infty$).

\medskip
It was conjectured in \cite{rc} that the coefficients of the
polynomial $(-1)^k F_\mu(\bmp;-\bmq)$ are \emph{nonnegative}, where
$-\bmq=(-q_1,\dots,-q_m)$. This conjecture was proved in \cite{rc} for
the case $m=1$, i.e., when $\lambda$ is a $p\times q$ rectangle,
denoted $\lambda=\ptq$.  For $w\in\sn$ let $\kappa(w)$ denote the
number of cycles of $w$ (in the disjoint cycle decomposition of
$w$). The main result of \cite{rc} was the following (stated slightly
differently but clearly equivalent).

\begin{theorem} \label{thm1}
Let $\mu\vdash k$ and fix a permutation
$w_\mu\in \fs_k$ of cycle type $\mu$. Then
  $$ F_\mu(p;q)=(-1)^k\sum_{uw_\mu=v}
      p^{\kappa(u)}(-q)^{\kappa(v)}, $$ 
where the sum ranges over all $k!$ pairs
$(u,v)\in\fs_k\times\fs_k$ satisfying $uw_\mu=v$.
\end{theorem}

To state our conjectured generalization of Theorem~\ref{thm1}, let
$\skm$ denote the set of permutations $u\in\fs_k$ whose cycles are
colored with $1,2,\dots,m$. More formally, if $C(u)$ denotes the set
of cycles of $u$, then an element of $\skm$ is a pair $(u,\varphi)$,
where $u\in\sk$ and $\varphi:C(u)\rightarrow [m]$. (We use the
standard notation $[m]=\{1,2,\dots,m\}$.) If
$\alpha=(u,\varphi)\in\skm$ and $v\in \sk$, then define a ``product''
$\alpha v=(w,\psi)\in\skm$ as follows. First let $w=uv$. Let
$\tau=(a_1,a_2,\dots,a_j)$ be a cycle of $w$, and let $\rho_i$ be the
cycle of $u$ containing $a_i$. Set
   $$ \psi(\tau)=\max\{
      \varphi(\rho_1),\dots,\varphi(\rho_j)\}. $$
For instance (multiplying permutations from left to right),
   $$ (\overbrace{1,2,3}^1)(\overbrace{4,5}^2)(\overbrace{6,7}^3)
          (\overbrace{8}^2)\cdot (1,7)(2,4,8,5)(3,5) =
     (\overbrace{1,4,2,6}^3)(\overbrace{3,7}^3)(\overbrace{5,8}^2). $$
Note that it a immediate consequence of the well-known formula
  $$ \sum_{w\in\fs_k} x^{\kappa(w)} = x(x+1)\cdots (x+k-1) $$
that $\#\skm = (k+m-1)_k$.  

\medskip
\textsc{Note.} The product $\alpha v$ does not seem to have
nice algebraic properties. In particular, it does not define an action
of $\sk$ on $\skm$, i.e., it is not necessarily true that $(\alpha u)v
= \alpha(uv)$. For instance (denoting a cycle colored 1 by leaving it
as it is, and a cycle colored 2 by an overbar), we have
   \beas [(\overline{1})(2)\cdot (1,2)]\cdot (1,2) & = &
       (\overline{1})(\overline{2})\\
    (\overline{1})(2)\cdot[(1,2)\cdot (1,2)] & = &
   (\overline{1})(2). \eeas 

\medskip
\indent Given $\alpha=(u,\varphi)\in\skm$, let
$\bm{p}^{\kappa(\alpha)}=
\prod_i p_i^{\kappa_i(\alpha)}$, where $\kappa_i(\alpha)$ denotes the
number of cycles of $u$ colored $i$, and similarly
$\bm{q}^{\kappa(\beta)}$, so
$(-\bm{q})^{\kappa(\beta)}= \prod_i (-q_i)^{\kappa_i(\beta)}$
We can now state our conjecture.

\begin{conjecture} \label{conj}
Let $\lambda$ be the partition of $n$ given by Figure~\ref{fig:nrec}.
Let $\mu\vdash k$ and fix a permutation $w_\mu\in \sk$ of cycle type
$\mu$. Then  
  $$ F_\mu(\bmp;\bmq)=(-1)^k\sum_{\alpha
      w_\mu=\beta} \bm{p}^{\kappa(\alpha)}(-\bm{q})^{\kappa(\beta)}, $$ 
where the sum ranges over all $(k+m-1)_k$ pairs
$(\alpha,\beta)\in\skm\times\skm$ satisfying $\alpha w_\mu=\beta$. 
\end{conjecture} 

\begin{example}
Let $m=2$ and $\mu=(2)$, so $w_\mu=(1,2)$. There are six pairs
$(\alpha,\beta)\in\sn^{(2)}$ for which $\alpha (1,2)=\beta$,
viz. (where as in the above Note an unmarked cycle is colored 1 and a
barred cycle 2), 
\medskip

\begin{center}
   \begin{tabular}{ccc} 
   $\alpha$ & $\beta$ & $\bm{p}^{\kappa(\alpha)} 
       \bm{q}^{\kappa(\beta)}$\\ \hline
    $(1)(2)$ & $(1,2)$ & $p_1^2q_1$\\
    $(\overline{1})(2)$ & $(\overline{1,2})$ & $p_1p_2q_2$\\
    $(1)(\overline{2})$ & $(\overline{1,2})$ & $p_1p_2q_2$\\
    $(\overline{1})(\overline{2})$ & $(\overline{1,2})$ & $p_2^2q_2$\\
    $(1,2)$ & $(1)(2)$ & $p_1q_1^2$\\
    $(\overline{1,2})$ & $(\overline{1})(\overline{2})$ & $p_2q_2^2$. 
   \end{tabular}
\end{center}
\medskip

\noindent It follows (since the conjecture is true in this case) that  
  $$  F_2(p_1,p_2;q_1,q_2) = -p_1^2q_1-2p_1p_2q_2
       -p_2^2q_2+p_1q_2^2+p_2q_2^2. $$
\end{example}

We can reduce Conjecture~\ref{conj} to the case $p_1=\cdots=p_m=1$;
i.e., $\lambda=(q_1,,q_2,\dots,q_m)$. Let 
  $$ G_\mu(\bm{p},\bm{q}) = \sum_{\alpha
      w_\mu=\beta} \bm{p}^{\kappa(\alpha)}\bm{q}^{\kappa(\beta)},
  $$  
so that Conjecture~\ref{conj} asserts that
$F_\mu(\bmp;\bmq) = (-1)^kG_\mu(\bm{p},-\bm{q})$.

\begin{proposition}
We have
  $$ \left.G_\mu(\bm{p},\bm{q})\right|_{q_{i+1}=q_i} =
   G_\mu(p_1,\dots,p_{i-1},p_i+p_{i+1},p_{i+2},\dots,p_m; $$
   \beq \hspace{1in}
   q_1,\dots,q_{i-1},q_i,q_{i+2},\dots,q_m). \label{eq:red} \eeq 
\end{proposition}

\proof 
Let $\alpha w_\mu=\beta$, where $\alpha,\beta\in\skm$ and $\mu\vdash
k$. If $\tau$ is a cycle of $\beta$ colored $i+1$ then change the
color to $i$, giving a new colored permutation $\beta'$. We can also
get the pair $(\alpha,\beta')$ by changing all the cycles in $\alpha$
colored $i+1$ to $i$, producing a new colored permutation $\alpha'$
for which $\alpha' w_\mu=\beta'$, and then changing back the colors of
the recolored cycles of $\alpha$ to $i+1$. Equation~(\ref{eq:red}) is
simply a restatement of this result in terms of generating functions.  
\qed

It is clear, on the other hand, that
  $$ \left.F_\mu(\bm{p},\bm{q})\right|_{q_{i+1}=q_i} =
     F_\mu(p_1,\dots,p_{i-1},p_i+p_{i+1},p_{i+2},\dots,p_m; $$
\vspace{-.35in}
   $$ \hspace{1in} q_1,\dots,q_{i-1},q_i,q_{i+2},\dots,q_m), $$
because the parameters $p_1,\dots,p_m; q_1,\dots,q_{i-1},q_i,q_i,
q_{i+2},\dots,q_m$ and
$p_1,\dots,p_{i-1},p_i+p_{i+1},p_{i+2},\dots,p_m;
q_1,\dots,q_{i-1},q_i,q_{i+2},\dots,q_m$ specify the same shape
$\lambda$. (Note that Proposition~\ref{prop:f} requires only $q_1\geq
q_2\geq \cdots \geq q_m$, not $q_1>q_2>\cdots>q_m$.)  Hence if
Conjecture~\ref{conj} is true when $p_1=\cdots=p_m=1$, then it is true
in general by iteration of equation~(\ref{eq:red}).

\medskip \textsc{Remarks.}  1. Conjecture~\ref{conj} has been proved by
Amarpreet Rattan \cite{rattan} for the terms of highest degree of
$F_k$, i.e., the terms of $F_k(\bmp; \bmq)$ of total degree $k+1$.
 
2. Kerov's character polynomials (e.g., \cite{g-r}) are
related to $F_k(\bm{p}; \bm{q})$ and are also conjectured to have
nonnegative (integral) coefficients. Is there a combinatorial
interpretation of the coefficients similar to that of
Conjecture~\ref{conj}?

\pagebreak

\end{document}